\documentclass[a4paper,12pt,oneside,reqno]{amsart}
\usepackage{amssymb,amsfonts,amsmath,amsthm,amscd, bbm}
\usepackage{graphicx, color}
\usepackage{mathrsfs}
\numberwithin{equation}{section}

\newcommand{\beq}{\begin{equation}}
\newcommand{\eeq}{\end{equation}}
\newcommand{\bea}{\begin{aligned}}
\newcommand{\eea}{\end{aligned}}
\newcommand{\bdm}{\begin{displaymath}}
\newcommand{\edm}{\end{displaymath}}
\newcommand{\barr}{\begin{array}}
\newcommand{\earr}{\end{array}}
\newcommand{\ben}{\begin{enumerate}}
\newcommand{\een}{\end{enumerate}}
\newcommand{\bde}{\begin{description}}
\newcommand{\ede}{\end{description}}

\newcommand{\R}{\mathbb{R}}

\newcommand{\N}{\mathbb{N}}

\newcommand{\PP}{\mathbb{P}}
\newcommand{\E}{{\mathbb{E}}}
\newcommand{\EE}{{\sf{E}}}

\newcommand{\defi}{\equiv} 

\newcommand{\be}{\beta}

\newcommand{\s}{\sigma}

\begin{document}

\title[From Parisi to Boltzmann]
{From Parisi to Boltzmann.}

\author[G. Kersting]{Goetz Kersting}
\address{Goetz Kersting \\ J.W. Goethe-Universit\"at Frankfurt, Germany.}
\email{kersting@math.uni-frankfurt.de}

\author[N. Kistler]{Nicola Kistler}
\address{Nicola Kistler \\ J.W. Goethe-Universit\"at Frankfurt, Germany.}
\email{kistler@math.uni-frankfurt.de}

\author[A. Schertzer]{Adrien Schertzer}
\address{adrien schertzer \\ J.W. Goethe-Universit\"at Frankfurt, Germany.}
\email{schertzer@math.uni-frankfurt.de}

\author[M. A. Schmidt]{Marius A. Schmidt}
\address{Marius A. Schmidt \\ University of Basel, Switzerland.}
\email{mariusalexander.schmidt@unibas.ch}

 \date{\today}

\thanks{ It is our pleasure to dedicate this work to Anton Bovier on the occasion of his $60^{th}$ birthday. }

\begin{abstract}
We sketch a new framework for the analysis of disordered systems, in particular mean field spin glasses, which is variational in nature and within the formalism of classical thermodynamics. For concreteness, only the Sherrington-Kirkpatrick model is considered here. For this we show how the Parisi solution (replica symmetric, or when replica symmetry is broken) emerges, in large but finite volumes, from a high temperature expansion to second order of the Gibbs potential with respect to order parameters encoding the law of the effective fields. In contrast with classical systems where convexity in the order parameters  is the default situation, the functionals employed here are, at infinite temperature, concave: this feature is eventually due to the Gaussian nature of the interaction and implies, in particular, that the canonical Boltzmann-Gibbs variational principles must be reversed. The considerations suggest that thermodynamical phase transitions are intimately related to the divergence of the infinite expansions.
\end{abstract}


\maketitle


\section{Introduction}
Despite the steady progress over the last decades in the rigorous treatment of mean field spin glasses, see \cite{talavol, panchy} and references therein, the Parisi solution \cite{mpv} still poses a number of deep conceptual questions. This is arguably due to the fact that the {\it physical content} of the theory still remains rather mysterious. In fact, the Parisi theory relies on reversed variational principles involving functional order parameters which are, from the point of view of classical thermodynamics, rather unorthodox. 

The goal of this paper is to shed some light on these puzzles: we shall introduce a Legendre formalism acting on finite volumes, and which rests on the insight that the order parameter of disordered systems must account for the randomness of the effective fields acting on the spins, and which survives the passage to the limit. As such, the order parameter must encode a {\it distribution}, and this in turn indeed suggest its "functional nature". Under this light, the replica symmetric solution of the SK-model, which is driven by a scalar order parameter, should be rather seen  as the exception to the rule, in presence of disorder; this is eventually due to the fact that in high temperature the limiting effective fields of the SK-model remain Gaussian: as such, their law is captured by few parameters only, the mean and the variance (in fact: only the variance). This statement must  sound like common place, especially when the SK is compared to, say, diluted models, where it is clear that no matter how large the temperature is, the order parameter must be given by a {\it function}, which indeed encodes the {\it law} of the effective fields see e.g. \cite[Chapter VI]{talavol} and references therein. 

The framework we wish to suggest finds its roots in a treatment of the Curie-Weiss model which is both extremely natural from a conceptual point of view, and  yet extremely challenging from the point of view of a rigorous implementation: it is based on high temperature expansions of the Legendre transformation of the (finite volume) free energy once a proper {\it Ansatz} is made for the effective magnetic fields acting on a spin due to the interactions. The key steps of the treatment are sketched in Section \ref{section_cw}. The consequences of such a point of view on disordered systems, such as the prototypical SK-model, are then worked out in Section \ref{section_disordered}, where both the replica symmetric phase and the 1RSB-phase are addressed. Anticipating the outcome of our considerations, the Parisi  solution, be it the replica symmetric- or the replica symmetry breaking version, turns out to correspond to the critical points of the Gibbs potential (Legendre transformation) expanded to second order in the inverse of temperature. We stress that the previous sentence concerns the Parisi solution, by this we understand the minimal values of the Parisi  functionals and not, as will become clear in the course of the discussion, the Parisi functionals themselves. The stability of the RS(B)-solution appears to be intimately related to the convergence of the Taylor series; in particular, it follows that thermodynamical phase transitions are due to  ("infrared") divergence of the high temperature expansions. 

Little/no emphasis on rigor is made in this paper: the goal is to address the question "what is the Parisi solution" within a transparent framework abiding to the principles of classical thermodynamics. A definite answer to this question shall naturally have consequences on a mathematical treatment, but this is left for future research. 

\section{Classical systems} \label{section_cw}
The Curie-Weiss (CW) model is an infinite-range system consisting of Ising spins $\boldsymbol \s = (\s_1, \dots, \s_N) \in \{\pm 1\}^N$ which interact through the Hamiltonian 
\beq 
H_N^\text{cw}( \boldsymbol \s) \defi  \frac{1}{N} \sum_{1\leq i < j \leq N} \s_i \s_j.
\eeq
The partition function to inverse temperature $\be = 1/kT>0$ (with $k$ the Boltzmann constant) and external field $h \in \R$, is given by
\beq 
Z_N^\text{cw}(\be,h) \defi E_o\left[ \exp\left( \be H_N^\text{cw}(\boldsymbol \s) + h \sum_{i\leq N} \s_i\right) \right],
\eeq 
where $E_o$ denotes average with respect to the uniform distribution $P_o(\boldsymbol \s) \defi 2^{-N}$. 

A key question concerns the infinite volume limit  of the free energy 
\beq 
F_N^\text{cw}(\be,h) \defi \frac{1}{N} \log Z_N^\text{cw}(\be,h). 
\eeq 
The free energy satisfies the fundamental principle 
\beq \bea \label{one}
& N F_N^\text{cw}(\be,h) = \max_{Q } \left\{ \int \left( \be H_N^\text{cw}(\boldsymbol \s) + h \sum_{i\leq N} \s_i \right) Q(d\s) - H(Q | P_o) \right\},
\eea \eeq
where $Q$ is any probability measure on $\{\pm 1\}^N$, the first term on the r.h.s. is the average energy with respect to $Q$, the {\it internal energy}, and the second term is the entropy of $Q$ relative to $P_o$.  The functional on the r.h.s. is, in finite volume, strictly concave in $Q$: the maximizer is unique, and given by the Gibbs measure 
\beq
\mathcal G_{N}^\text{cw}(\boldsymbol\s) \defi \exp(\be H_N^\text{cw}(\boldsymbol \s) + h \sum_{i\leq N} \s_i)\Big/Z_N(\be,h).
\eeq 

\subsection{Gibbs potential, and high temperature expansions}
A number of elementary approaches have been developed in order to compute the limiting free energy of the CW-model, see \cite{bovier} and references therein. Here we shall discuss informally an approach via  high temperature expansions, much akin to the setting of Plefka \cite{plefka} for the TAP-analysis \cite{tap} of the SK-model. The 'spirit' of the approach shall motivate and justify the treatment of disordered systems which is  sketched in Section \ref{section_disordered}.

Let us shorten $\mathcal H_N^\text{cw}(\s) \defi \be H_N^\text{cw}(\s)$. For $\alpha \in \R$, and $\boldsymbol \varphi= \{\varphi_i\}_{i =1 \dots N} \in \R^N$, we introduce the (normalized) functional $G_N(\alpha, \boldsymbol \varphi)$ according to 
\beq \label{g}
N G_N(\alpha, \boldsymbol \varphi) = \log E_o\left[ \exp\left(\alpha \mathcal H_N^\text{cw}(\s)+ h \sum_{i\leq N} \s_i + \sum_{i\leq N} \varphi_i \s_i\right) \right]
\eeq
Remark that by Jensen's inequality, the map  $\boldsymbol \varphi \mapsto G_N(\alpha, \boldsymbol \varphi)$  is, in fact, {\it convex}. In particular, the Legendre transformation is well defined:
\beq \label{vp}
G_N(\alpha, \boldsymbol m)^\star \defi \sup_{\boldsymbol \varphi \in \R^N}  \sum_{i\leq N} \varphi_i m_i - G_N(\alpha, \boldsymbol \varphi)\,.
\eeq
The functional $G_N^\star(\alpha, \boldsymbol m)$ is typically referred to as {\it Gibbs potential} (or {\it Helmoltz free energy}: here and below, we shall adopt the former terminology). Again by convexity, the ${}^\star$-operation is an involution, i.e. with the property that $G_N = \left( G_N^\star \right)^\star$. Since by construction $G_N(1, \boldsymbol 0)$ coincides with the free energy, we therefore have that 
\beq \bea
F_N(\be, h) = \sup_{\boldsymbol m \in \R^N}  \left\{ -G_N( 1, \boldsymbol m)^\star \right\} \,.
\eea \eeq
The thermodynamic variables $\boldsymbol m \in \R^N$ are dual to the magnetic fields $\boldsymbol \varphi$, and correspond, upon closer inspection, to the magnetization: indeed, denoting by $\left< \right>_\alpha$ Gibbs measure with respect to the Hamiltonian appearing in \eqref{g}, one immediately checks by solving the variational principle \eqref{vp} that the fundamental relation holds 
\beq \label{mm}
\left< \s_i \right>_\alpha = m_i \,.
\eeq
(In particular, we see from the above that $\boldsymbol m \in [-1,1]^N$).  The idea is to now proceed by Taylor expansion of the Gibbs potential, 
\beq
- G_N(\alpha, \boldsymbol m)^\star = \sum_{k=0}^\infty  \frac{d^k}{d\alpha^k}\Big( - G_N(\alpha, \boldsymbol m)^\star \Big) \Big|_{\alpha=0} \frac{\alpha^k}{k!}\,,
\eeq
and to evaluate this in $\alpha=1$. The calculation of the Taylor-coefficients considerably simplifies in $\alpha=0$, as one only needs to compute "spin-correlations "under the non-interacting Hamiltonian $\sum_{i\leq N} \varphi_i \s_i$. One immediately checks that the $0^{th}$-term of the expansion is given by 
\beq 
- G_N(0, \boldsymbol m) = \frac{1}{N} \left( h \sum_{i\leq N} m_i - \sum_{i\leq N} I(m_i) \right)
\eeq
with 
\[
I(x) \defi \frac{1+x}{2}\log(1+x)+ \frac{1-x}{2}\log(1-x), \qquad x \in [-1,1],
\]
the rate function for Ising spin, i.e. the entropic cost for fixing the spin-magnetizations to the prescribed values \eqref{mm}. 

The first derivative in $\alpha=0$ is also elementary: it gives a contribution 
\beq \label{fd}
- G_N'(0, \boldsymbol m) = \frac{\be}{N^2} \sum_{i<j}^N m_i m_j,
\eeq
which we immediately recognize as the internal energy under the non-interacting Hamiltonian.

The  higher order derivatives all give a contribution which is irrelevant in the $N\to \infty$ limit, {\it provided} that  the series expansion is absolutely convergent up to $\alpha=1$. According to \cite{plefka}, this is the case for $\boldsymbol m$ satisfying the restriction of the mean field theory for an Ising ferromagnet, to wit:
\beq \label{ht}
\be\left(1-\frac{1}{N} \sum_{i=1}^N m_i^2 \right) < 1 .
\eeq

All in all, the "high temperature expansion" of the Gibbs potential with respect to the magnetization as order parameter, leads to the expression
\beq \bea \label{vp}
F_N^\text{cw}(\be, h) &= \hat \sup_{\boldsymbol m} \left\{  \frac{h}{N} \sum_{i\leq N} m_i + \frac{\be}{N^2} \sum_{i<j}^N m_i m_j - \frac{1}{N} \sum_{i\leq N} I(m_i)  + O(1/N) \right\},
\eea \eeq
where the supremum is over magnetizations satisfying \eqref{ht}.

It is of course a simple task to solve the above variational principle: by symmetry, one expects the supremum to be achieved in $m_i= m \in [-1,1]$ for $i=1\dots N$, in which case one gets
\beq \label{lim_f}
\lim_{N\to \infty} F_N(\be, h) = \max_{\be(1-m^2)< 1} h m +\frac{\be}{2} m^2 - I(m),
\eeq
which is the well-known solution of the CW, see e.g. \cite{bovier}. The maximization, and the concavity properties of the limiting $m$-functional on the r.h.s. above, are of course in complete agreement with the finite volume variational principle \eqref{one}. 

One can hardly overstate that \eqref{lim_f}, although emerging from a high temperature expansion, is valid for {\it any} $\be$: this is due to the mean field character of the CW-model, and the fact that in low temperature, the Gibbs measure concentrates on finitely many pure states (in fact only one pure state if $h\neq 0$, and  two distinct pure states detecting the symmetry breaking under spin-flips if $h=0$) which are all effectively in high temperature. 

\section{Disordered Systems} \label{section_disordered}
Archetypical mean field spin glasses are the socalled $p$-spin models (and mixtures thereof). Here, we will stick to the $p=2$ case, also known as the Sherrington-Kirkpatrick model \cite{sk}. Consider to this end the space of Ising configurations in finite volume $N$, i.e. $\Sigma_N \defi \{\pm 1\}^N$, the space of Ising configurations. Consider also 
independent standard Gaussians $\{ g_{ij}\}_{1 \leq i < j \leq \N}$ issued on some probability space $(\Omega, \mathcal F, \PP)$ and the Hamiltonian 
\beq 
H_N^\text{sk}(\s) \defi \frac{1}{\sqrt{N}} \sum_{i<j} g_{ij} \s_i \s_j,
\eeq
We define the finite volume quenched free energy in magnetic field $h\in \R$ according to 
\beq
F_N(\be, h) \defi \frac{1}{N} \E \log E_o\left[ \exp\left( \be H_N^\text{sk}(\s)+h \sum_{i\leq N} \s_i \right) \right]\,.
\eeq
In line with the previous section, we would like to find an order parameter, and a Legendre structure behind the free energy on which we can perform high temperature expansions.  Contrary to the CW-model, this is no simple task. We shall begin in the following section with the simplest case: the so-called replica symmetry [RS] phase. The more involved case of replica symmetry breaking [RSB] will be addressed in Section \ref{rsb_sec} below. 

\subsection{High temperature, or: the RS-Legendre transformation} \label{rs}
Before addressing the SK-model, we shall ponder on the Legendre transformation behind the Curie-Weiss model, addressing the question "what is the operation really doing". The key insight here is that, due to the mean field character of the interaction, the Gibbs measure should split, in the infinite volume limit, into a convex combination of pure states which are all at high temperature.
In other words the Gibbs measure should be a mixture of ergodic components/ product measures. As it turns out, an approximate product measure structure already hides behind the finite volume CW-Hamiltonian. To see this, let us denote by $\left< \right>_{\be, h, N}$ the expectation under the CW-Gibbs measure, and introduce the 
{\it effective magnetic fields} 
\beq 
h_i(\s) \defi h+ \frac{\be}{N} \sum_{j \neq i}^N \s_j.
\eeq  
With these notations we can re-write the CW-Hamiltonian as 
\beq
\be H_N^{\text{cw}}(\s) + h \sum_{i\leq N} \s_i = \sum_{i=1}^N h_i(\s) \s_i
\eeq
Naturally (yet perhaps with a grain of hindsight), we expect the free energy to be carried by configurations for which 
\beq
h_i(\s) \approx \left< h_i(\s) \right>_{\be, h, N}^{(i)}, 
\eeq
where $\left< \right>_{\be, h, N}^{(i)}$ denotes Gibbs-expectation with respect to a CW-Hamiltonian with a cavity in the site $i$. Denoting by $m_j^{(i)} \defi \left< \s_j \right>_{\be, h, N}^{(i)}$ and
$m^{(i)} \defi \frac{1}{N} \sum_{j \neq i} m_j^{(i)}$ we would thus have the "approximation"
\beq
\be H_N^{\text{cw}}(\s) + h \sum_{i\leq N} \s_i \approx \sum_{i\leq N} \left(  \be m^{(i)} + h \right) \s_i
\eeq
Replacing now $ \be m^{(i)} +h =: \varphi_i$, we thus see the emergence of the Hamiltonian $\sum_i \varphi_i \s_i$ as in \eqref{g}: the (double) Legendre transformation and the high temperature expansions thus identify self-consistently the effective fields $\varphi's$ leading to the correct mean magnetizations of the spins of the CW-model. Since any product measure on Ising spins is evidently uniquely characterized by the magnetization, the introduction of the $\varphi$-fields as order parameters is thus justified in rather natural terms. 

The above take also applies to disordered Hamiltonians such as the SK-model. Again due to the mean field character, one expects the quenched Gibbs measure of the SK-model to decompose in the limit into a mixture of product measures. (Contrary to the ordered case, however, one naturally expects the coefficients of the mixture to be also {\it random}). We now proceed through a sequence of radical assumptions, which one may expect to hold for small $\be$ (the regime of low correlations), and which are at the very basis of any treatment in the physical literature (see \cite{mpv} and references therein): 
\begin{itemize}
\item[i)] only one pure state contributes to the Gibbs measure. Under this assumption, a re-run of the above considerations leads to assume that the SK-Hamiltonian will be safely "approximated" by 
\beq
\be H_N^{\text{sk}}(\s)+ h \sum_{i =1}^N \s_i \approx \sum_{i=1}^N \left( \frac{\be}{\sqrt{N}} \sum_{j\neq i} g_{ij} m_j^{(i)} + h \right) \s_i
\eeq
\end{itemize}
The situation is however more delicate than in the CW-model: the scaling of the effective field is of order $\sqrt{N}$, in which case no strong concentration in the limit $N\to \infty$ can be expected. At the very best, only weak limits (in the form of a central limit theorem) can be expected to kick in, and we shall indeed assume this to be the case. Precisely:
\begin{itemize}
\item[ii)] We put forward the working assumption that 
\beq
\lim_{N\to \infty} \frac{\be}{\sqrt{N}} \sum_{j \neq i} g_{ij} m_j^{(i)} = \be \sqrt{q} g_i, 
\eeq
{\it weakly}, where $g_i \sim \mathcal N(0,1)$ is a standard Gaussian, and $$q \defi \lim_{N\to \infty} \frac{1}{N} \sum_{j \neq i} \left( m_j^{(i)} \right)^2 \in [0,1],$$ assuming of course the limit exists. In line with the CW-model, we may here expect the {\it variance} of the effective fields to play the role of order parameter. 
\item[iii)] We furthermore assume the simplest possible covariance structure among the $g's$, to wit: we shall assume these to be {\it independent} from one another. 
\end{itemize}
All in all, items {i-iii)} suggest to introduce the following functional: 
\beq \bea \label{1map}
q \in \R_+ \mapsto \Phi_\alpha(q) & \defi  \frac{1}{N} \E \log E_o \exp\left( \sqrt{\alpha} \mathcal H_N(\s) + h\sum_{i\leq N} \s_i+ \be \sqrt{q} \sum_{i\leq N} g_i \s_i \right)\,,
\eea \eeq
where $\mathcal H_N(\s) \defi \be H_N^{\text{sk}}(\s)$. This functional plays the essential role  of \eqref{g}, once the disordered nature of the Hamiltonian is taken into account. (The reader will notice that contrary to \eqref{g}, we have introduced here a square root dependence on the $\alpha$ parameter, in line with the square root dependence on $q$).  We will refer to \eqref{1map} as the {\it RS-Legendre functional}. 

Something curious is happening: for $\alpha=0$, the map $q \mapsto \Phi_0(q)$ is {\it not} convex, but {\it concave}. It is not true that concavity holds for all $\alpha$, but this shouldn't bother us: the goal being to shed light on the Parisi solution, it suffices to consider the {\it concave} Legendre transform 
\beq
\Phi_\alpha^\star(q^\star) = \min_{q \in [0,1]} q q_\star - \Phi_\alpha(q),
\eeq
with $q^\star$ being the dual of the parameter $q$, and appeal to the fact that double Legendre leads to the {\it concave envelope}, which upper bounds the original function:
\beq
\Phi_\alpha(q) \leq \left(\Phi_\alpha^\star\right)^\star(q) = \min_{q^\star} q q^\star - \Phi_\alpha^\star\left(q^\star\right)
\eeq
In particular we see that the {finite volume} free energy is upper  bounded by
\beq \label{sk_prin}
F_N(\be, h) = \Phi_1(0) \leq \left(\Phi_1^\star\right)^\star(0) = \min_{q^\star} \max_{q} \Phi_1(q) -qq^\star.
\eeq
As in the case of the CW-model,   the idea is now to Taylor-expand the inner maximization, which we refer to as the "replica symmetric Gibbs potential", around $\alpha=0$, and to evaluate this in $\alpha=1$. Shortening 
\beq
\tilde \phi(\alpha, q^\star) \defi \max_{q} \Phi_\alpha(q) -qq^\star,
\eeq 
we thus seek the expansion
\beq
\tilde \phi(1, q^\star) = \tilde \phi(0, q^\star)+ \tilde \phi'(0, q^\star)+ \frac{1}{2} \tilde \phi''(0,q^\star)+ \dots
\eeq
Towards this goal, we shall slightly deviate from the CW-approach insofar we make use of the 
$q's$ as thermodynamical variables, the reason being that no explicit expression of the expansion can be reached in the $q_\star$ formulation, and this would render the analysis cumbersome. By a slight abuse of notation, we will write $\tilde \phi(\alpha, q) \defi \tilde \phi(\alpha, q^\star(\alpha, q))$, with $q^\star(\alpha, q)$ given by 
\beq
q^\star(\alpha, q) \defi \partial_q \Phi_\alpha(q). 
\eeq
In $\alpha=0$, and by Gaussian P.I., this takes the form,   
\beq \label{non_lin}
q^\star(0, q) = \frac{\be^2}{2} \left\{ 1- \EE \tanh^2(h+\be \sqrt{q} g)\right\},
\eeq
where $g$ is a standard Gaussian. Therefore, 
\beq \bea \label{0}
\tilde \phi(0, q) &= \Phi_0(q) - q \cdot q^\star(0, q) \\
&=  \EE \log \cosh(h+\be \sqrt{q} g) - \frac{\be^2}{2} q \left\{ 1- \EE \tanh^2(h+\be \sqrt{q} g)\right\}\,.
\eea \eeq
The first $\alpha$-derivative is steadily computed: 
\beq \bea \label{first_gip}
\frac{d }{d\alpha} \tilde \phi(\alpha, q) &= \frac{\partial}{\partial_\alpha}  \tilde \phi(\alpha, q)  \hspace{4.3cm} \text{(extremality of $q$)} \\
&= \frac{1}{ 2 N \sqrt{\alpha}} \E \left< \mathcal H_N(\s) \right>_\alpha \\
& = \frac{\be}{ 2 N \sqrt{\alpha N}} \sum_{i< j}\E g_{ij} \left< \s_i \s_j \right>_\alpha \\
& = \frac{\be^2}{2 N^2} \sum_{i< j}\left( 1 - \E \left< \s_i \s_j \tau_i \tau_j \right>^{\otimes 2}_\alpha \right)   \hspace{1.2cm} \text{(Gaussian P.I.)}
\eea \eeq
where $\left< \right>_\alpha^{\otimes 2}$ stands for quenched Gibbs measure with respect to the Hamiltonian \eqref{1map} over the replicated configuration space, and  $(\s, \tau) \in  \Sigma_N \times \Sigma_N$.
In $\alpha=0$, and thanks to the decoupling of the spins, this expression can be explicitely computed: the upshot is 
\beq \bea \label{1}
\frac{d }{d\alpha} \tilde \phi(\alpha, q) \Big|_{\alpha=0} & = \frac{\be^2}{2N^2} \sum_{i< j}\left\{ 1-  \EE[\tanh^2(h+\be \sqrt{q} g)]^2 \right\} \\
& = \frac{\be^2}{4}\left\{ 1-  \EE[\tanh^2(h+\be \sqrt{q} g)]^2 \right\} + o(1) \quad (N\to \infty).
\eea \eeq
To lighten notation, we shorten
\beq \bea \label{exp_gibbs}
\text{G}_{\text{rs}}(q) & \defi \tilde \phi(0, q) + \tilde \phi'(0, q) \stackrel{\eqref{0}, \eqref{1}}{=} \EE \log \cosh(h+\be \sqrt{q} g)  + \\
& \quad + \frac{\be^2}{4}\left\{ 1-  \EE[\tanh^2(h+\be \sqrt{q} g)] \right\} \left\{ 1+  \EE[\tanh^2(h+\be \sqrt{q} g)] -2 q \right\}\,,
\eea \eeq
for the first order expansion of the RS-Gibbs potential. Altogether, 
\beq  \label{final_al}
\lim_N F_N(\be, h)  \leq \min_q \text{G}_{\text{rs}}(q) +  \text{(higher derivatives)}.
\eeq
The higher order derivatives in \eqref{final_al} correspond to  spin-correlations (and moments thereof) under the non-interacting Hamiltonian ($\alpha=0$); as such they can be explicitely computed.  One checks that  the contribution of the second derivative already drops down to $O(1/N)$. Unfortunately, determining the radius of convergence of the infinite expansions turns out to be enormously challenging, but in analogy with the CW-model,  we may expect that for some $(\beta, h)$-regime the infinite series are absolutely convergent, in which case the contribution of  the  higher orders would be irrelevant in the thermodynamical limit, and we would get
\beq \label{final}
\lim_{N} F_N(\be, h)  \leq \min_q \text{G}_{\text{rs}}(q)\,.
\eeq
The variational principle on the r.h.s. of \eqref{final} can be easily solved: the minimum is attained in the solution of the (well-known) equation
\beq \label{fixpoint}
 q= \EE\tanh^2(h+\be \sqrt{ q}g),
\eeq
in which case, denoting by $\hat q$ the solution\footnote{It is known that the solution is unique for small $\beta$, see e.g. \cite{panchenko} and references therein.} of the fixpoint of \eqref{fixpoint}, we obtain
\beq \label{free}
\lim_N F_N(\be, h) \leq \EE \log \cosh(h+\be \sqrt{\hat q} g) +\frac{\be^2}{4}\left( 1- \hat q \right)^2.
\eeq 
The r.h.s. is the celebrated RS-solution of the SK-model \cite{sk, mpv} which we know is correct for small $\beta$, see \cite{guerra_toninelli, tala, panchenko} and references therein. Here is a first bottom line: 

\begin{center}
{\sf the RS Parisi solution  is a first-order expansion of the \\
RS-Gibbs potential w.r.t. to the variance of the effective fields}. 
\end{center} 
${}$

The above wording, in particular the "first order" characterization,  is in line with the derivation but things should nevertheless be taken with a grain of salt. Indeed, the mean field correction as given in the second line of \eqref{exp_gibbs} is the outcome of an additional Gaussian partial integration \eqref{first_gip} which, in particular, increases the order of  derivatives by one. In other words, what seems to be a first order is, in fact, a second order, and this in agreement
with the fact that the mean field correction is {\it quadratic} in $\beta$. This stands of course in contrast with classical systems such as the CW-model, cfr. \eqref{lim_f} and in particular the linear $\beta$-dependence of the mean field correction in that case. The SK-quadratic correction is, in fact, the {\it signature} of disorder. 

Somewhat unexpectedly, the Legendre/high temperature expansions lead to variational principles which, to our knowledge, have not appeared in the literature so far. In fact, \eqref{final} does {\it not} coincide with the canonical formulation for the RS-phase of the SK-model. This observation leads to some intriguing insights concerning the high temperature phase, which we shall briefly discuss.  Let us shorten
\beq
F_\text{rs}(q) \defi \EE \log \cosh(h+\be \sqrt{q} g) +\frac{\be^2}{4}\left( 1- q \right)^2, 
\eeq
for the "classical" RS-functional of the SK-model.  It is then known \cite{guerra, tala} that 
\beq \label{boh}
\lim_N F_N(\be, h) = \inf_{q} F_\text{rs}(q)\,,
\eeq
provided $\beta$ is small enough. The functionals $G_{\text{rs}}$ and $F_\text{rs}$ are quite similar, and yet manifestly different: they achieve the {\it same} minimal value for $q= \hat q$ (the solution of the fixed point equation \eqref{fixpoint}) but surprisingly, the expansion of the Gibbs potential typically lies lower than the RS-functional, i.e.
\beq
G_{\text{rs}}(\hat q) = F_\text{rs}(\hat q)  
\eeq
for $\hat q$ solution of \eqref{fixpoint}, but 
\beq
G_{\text{rs}}(q) < F_\text{rs}(q), 
\eeq
strictly, for $q \neq \hat q$. Below are some numerical plots depicting the state of the matter: in each case, the  parameters $\beta, h$ are chosen to lie {\it below} the AT-line, hence the system is in the alleged high temperature/replica symmetric phase.

\begin{figure}[h]
    \centering
    \includegraphics[width=0.3\textwidth]{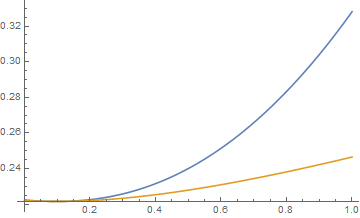} $\hspace{1cm}$ 
    \includegraphics[width=0.3\textwidth]{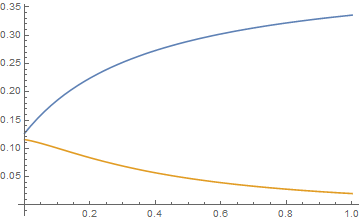}
    \caption{On the left: $F_\text{rs}$-functional (blue), $G_{\text{rs}}$-functional (yellow), for $\be=0.9$ and $h=0.2$. Both functions are manifestly convex, as confirmed by the figure on the r.h.s. which plots the corresponding second derivatives: both are positive for $q\in [0,1]$.}
    \label{fig:1}
\end{figure}

\begin{figure}[h]
    \centering
    \includegraphics[width=0.3\textwidth]{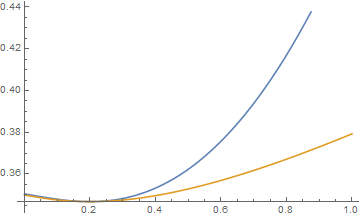} $\hspace{1cm}$ 
    \includegraphics[width=0.3\textwidth]{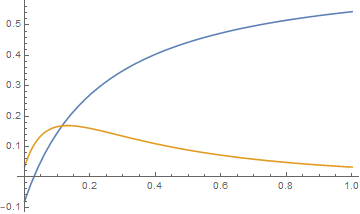}
    \caption{Same as above but $\be=1.15$ and $h=0.2$. The second derivative of the Gibbs potential (yellow) is always positive: the functional is thus globally convex. This is however not the case for the $F_\text{rs}$-functional: for small $q$-values the second derivative becomes negative, detecting a loss of convexity. }
    \label{fig:2}
\end{figure}

\begin{figure}[h]
    \centering
    \includegraphics[width=0.3\textwidth]{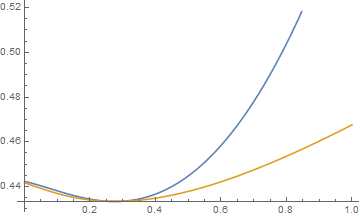} $\hspace{1cm}$ 
    \includegraphics[width=0.3\textwidth]{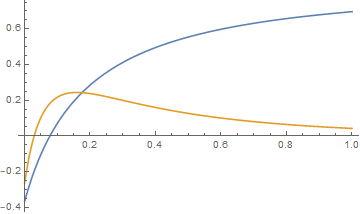}
    \caption{Same as above but $\be=1.3$ and $h=0.2$. Remark that both second derivatives become negative for small $q$-values: for such parameter-values, both functionals thus lose (global) convexity. }
    \label{fig:3}
\end{figure}

These plots help to visualise some underlying conceptual issues. First, concerning the (expansion of the) Gibbs potential: by its very derivation, whenever global convexity is lost (such is the case, e.g, in Figure \ref{fig:3} above) one must unequivocally deduce that the infinite expansions have stopped to converge globally. Dropping such contributions is therefore unjustified, and potential source of inconsistencies. (Here we are only rephrasing Plefka's deep insight \cite{plefka} that convexity properties of the limiting functionals are necessary, yet not necessarily sufficient conditions for the convergence of infinite expansions). 

Concerning the $F_{\text{rs}}$-functional, the conclusions are more vague, but one can safely say that whenever the functional is no longer convex (such is the case in Figure \ref{fig:2} and \ref{fig:3}) one must deduce that the {\it bona fide} variational principle \eqref{boh} can hardly be the outcome of a sound thermodynamical principle. This, of course, does {\it not} imply that the minimal value yields a wrong free energy\footnote{We shall mention in passing that these considerations might also have some consequences on the delicate, and to these days still debated validity of the AT-line \cite{at}: as the plots show, there are $(\beta, h)$-regimes below the AT-line where low temperature behavior is already "hiding" behind the RS-solution. One is thus presumably better off by considering the $K$-RSB formulation \cite{mpv} (potentially with $K=\infty$), which can still collapse to the RS-value. Indeed, it is known \cite{ac_two} that the full Parisi variational principle is convex in Parisi's functional order parameter, the latter belonging to the convex set of increasing functions on the unit interval.  Since the local minimum of a globally convex functional is also the global minimum, to (dis)prove the validity of the AT-line amounts to (dis)proving that the RS-solution is a local minimum.}. The crux of the matter is that contrary to \eqref{final}, which emerges from an expansion of the Gibbs potential, it is  unclear (to us) what the content of the variational principle \eqref{boh} is. \\

\subsection{Low temperature, or: the 1RSB-Legendre transformation} \label{rsb_sec}
It can hardly be stated too strongly that the Legendre transformations employed in the high temperature regime are nothing but an {\it Ansatz}: this is arguably self-evident in the SK-case \eqref{1map}, but it is nonetheless true for the "trivial" CW-model \eqref{g}. According to the above discussion, a fundamental insight in the case of the SK-model is a judicious guess concerning the {\it law} of the effective fields. In the high temperature phase, we have made the simplest possible choice, which we may loosely summarize in "what is Gaussian remains Gaussian". This, however, can hardly be correct for any $\beta, h$: at the very latest when the map $q \mapsto G_{\text{rs}}(q)$ ceases to be concave (as we have seen this definitely  happens for $\beta$ close to the AT-line) we must deduce that the expansions have already stopped to converge absolutely, in which case the treatment comes to a stall. In order to find a  suitable replacement of the Legendre functional \eqref{1map}, one seeks a better {\it Ansatz} for the law of the effective fields which also takes into account the existence of many pure states. This is of course a difficult issue, and we have no simple/easy-to-grasp justification for the choice we will make below and which is, not surprisingly, in line with the Parisi solution \cite{mpv, guerra, ass}. As explained in \cite[Chapter V.2]{mpv}, and due to the presence of many pure states in the low temperature regime, the law of the effective fields is given by certain "tilted Gaussian" measures. To take this into account we consider parameters $\boldsymbol m = (m_0, m_1)$ and $\boldsymbol q = (q_0, q_1)$ where 
$0\leq m_0 \leq m_1 \leq 1$, $0\leq q_0 \leq q_1 \leq 1$, and introduce independent standard Gaussians $g^{(0)}_{i}, g^{(1)}_{i}, i=1\dots N$. Consider then the following functional  
\beq \bea \label{2map}
& \Phi_{\alpha}(\boldsymbol q; \boldsymbol m) \defi  \\
& \frac{1}{N} \E\left[ \frac{1}{m_0} \log \EE_0\left\{ \exp \frac{m_0}{m_1} \log \EE_1 \left[ E_o\left[ e^{ \sqrt{\alpha} \mathcal H_N(\s)+ \sum_{i\leq N}\left(h+\be \sqrt{q_0} g_i^{(0)} + \be \sqrt{q_1-q_0} g_i^{(1)}\right] \s_i } \right]^{m_1}
\right] \right\} \right]\,.
\eea \eeq
In the above: $\E$ denotes expectation of the SK-disorder; $\EE_0$ and $\EE_1$ denote 
$\{g_i^{(0)} \}_{i\leq N}$ and respectively $\{ g^{(1)}_i \}_{i\leq N}$ expectations. Finally, as before, $E_o$ denotes expectation with respect to coin tossing measure on $\Sigma_N$.

We point out that a very similar object has been first employed by Guerra in the groundbreaking \cite{guerra}.  

The functional \eqref{2map} should be seen as a generalization of \eqref{1map}: one checks that 
\beq
\lim_{m_0\downarrow 0} \Phi_{\alpha}(q_0=q,q_1=q; \boldsymbol m) = \Phi_\alpha(q).
\eeq
By analogy with \eqref{1map}, we thus refer to \eqref{2map} as the {\it 1RSB Legendre functional}. 

Remark that for $\alpha=1$, 
\beq
F_N(\be, h) = \Phi_{1}(\boldsymbol 0; \boldsymbol m) \,,
\eeq
and that in $\alpha=0$ the mapping is concave in $q_0, q_1$. This again suggests the use of concave Legendre transformations. We shall however perform {\it no} Legendre transformation on $\boldsymbol m$: these parameters should be thought as being given, and fixed\footnote{The (probably deep) reason eludes us: quite simply, regarding the $\boldsymbol m's$ as fixed parameters leads seamlessly, as we are going to see, to the Parisi 1RSB solution; on the other hand, performing a Legendre transformation also on these parameters yields unwieldy/unrecognizable variational principles, see also the discussion in the last Section.}. Precisely, we consider the concave Legendre transformation 
\beq
\Phi_{\alpha}^\star(q_0^\star, q_1^\star; \boldsymbol m) \defi \min_{q_0^\star, q_1^\star} \left\{ q_0 q_0^\star+q_1 q_1^\star - \Phi_{\alpha}(\boldsymbol q; \boldsymbol m) \right\}
\eeq
Since double Legendre leads to the concave hull, we thus have 
\beq
F_N(\be, h) \leq \left( \Phi_{1}^\star\right)^\star(0, 0; \boldsymbol m) = \min_{q_0^\star, q_1^\star} \left\{ - \Phi_{1}^\star(q_0^\star, q_1^\star; \boldsymbol m) \right\}. 
\eeq
As in the replica symmetric case, the idea is now to Taylor expand the Gibbs potential around $\alpha=0$, and evaluate it in $\alpha=1$. The analysis is elementary yet somewhat cumbersome due to the higher complexity of the involved functions. Anticipating the upshot, we will see that stopping the expansion to first order yields a functional which, in a specific choice of the $(\boldsymbol q, \boldsymbol m)$-variables, collapses to the 1RSB Parisi solution. 

The first step is thus to identify 
\beq \label{varprin_2}
\arg \min_{q_0^\star, q_1^\star} \left\{ q_0 q_0^\star+q_1 q_1^\star - \Phi_{0}(q_0, q_1; \boldsymbol m)\right\}\,,
\eeq
and for this some notation is needed. We define the functions
\beq
f(g_0, g_1; \boldsymbol q) \defi \log \cosh(h+\be \sqrt{q_0} g_0 + \be \sqrt{q_1-q_0} g_1)
\eeq
(we omit $\be, h$, which are fixed, to lighten notation), and  
\beq
f_0(g_0; \boldsymbol q, \boldsymbol m) \defi \frac{1}{m_1} \log \EE_1\left[\exp m_1 f(g_0, g_1; \boldsymbol q)\right]
\eeq
where $\EE_1$ stands for integration w.r.t. $\mathcal N(dg_1)$, the standard Gaussian. We also introduce the following densities:  
\beq
\nu_1(g_0, dg_1) \defi \frac{\exp m_1 f(g_0, g_1; \boldsymbol q, \boldsymbol m)}{\EE_1 \exp m_1 f(g_0, g_1;  \boldsymbol q, \boldsymbol m)} \mathcal N(dg_1), 
\eeq
\beq
\nu_0(dg_0) \defi \frac{\exp m_0 f_0(g_0, g_1; \boldsymbol q, \boldsymbol m)}{\EE_0 \exp m_0 f_0(g_0; \boldsymbol q, \boldsymbol m)} \mathcal N(dg_0).
\eeq
We omit the dependence on $\boldsymbol q, \boldsymbol m$ to lighten otherwise heavy notation, and remark for later use that
\beq \label{lim_m}
\lim_{m_0 \downarrow 0} \nu_0(dg_0) = \mathcal N(dg_0),
\eeq
i.e. the density collapses to the standard Gaussian. 
We then set $\nu \defi \nu_0 \otimes \nu_1$, which is a density on $\R^2$. Given a function $G: \R^2 \to \R, (g_0, g_1) \mapsto G(g_0, g_1)$ we shall write 
\beq
\nu\left[ G(g_0, g_1) \right]
\defi \int \left(\int G(g_0, g_1) \nu_1[g_0, dg_1] \right) \nu_0(dg_0)\,.
\eeq
Finally, we shorten
\beq
\text{Tanh}(g_0, g_1) \defi \tanh\left( h + \be \sqrt{q_0} g_0 + \be \sqrt{q_1-q_0} g_1 \right)\,.
\eeq
With these notations we thus have that 
\beq \bea
\Phi_{0}(\boldsymbol q; \boldsymbol m) &= \frac{1}{m_0} \log \EE_0\left[ \exp\left( \frac{m_0}{m_1} \log \EE_1\left[ \exp m_1 f(g_0, g_1; \boldsymbol q, \boldsymbol m)\right] \right) \right] \\
&= \frac{1}{m_0} \log \EE_0\left[ \exp\left( m_0 f_0(g_0; \boldsymbol q, \boldsymbol m)\right) \right]
\eea \eeq
Coming back to the variational principle \eqref{varprin_2}, taking the $q_0$ derivative and by Gaussian P.I. we see that $q_0$ and $q_1$ must solve, in $\alpha=0$, the equations
\beq \label{q0}
q_0^\star = \frac{\be^2}{2}(m_0-m_1)\nu_0 \left[\nu_1\left[ \text{Tanh}(g_0, g_1) \right]^2\right]\,,
\eeq
\beq \label{q1}
q_1^\star = \frac{\be^2}{2}-\frac{\be^2}{2}(1-m_1) \nu\left[  \text{Tanh}^2(g_0, g_1) \right]\,.
\eeq
As in the RS-case, these equations cannot be solved explicitely: we shall therefore use $q_0, q_1$ as  thermodynamical variables. We thus replace $q_0^\star = q_0^\star(q_0, q_1) , q_1^\star = q_1^\star(q_0, q_1)$ with their representations \eqref{q0}, \eqref{q1}, and by slight abuse of notation, we set 
\beq 
\Phi_0^\star(q_0, q_1, \boldsymbol m) \defi \Phi_{0}^\star(q_0^\star(q_0, q_1), q_1^\star(q_0, q_1); \boldsymbol m).
\eeq
It then holds 
\beq \bea
\Phi_0^\star(q_0, q_1) & = \frac{\be^2}{2}q_0(m_0-m_1)\nu_0\left[ \nu_1\left[ \text{Tanh}(g_0, g_1) \right]^2\right] +\frac{\be^2}{2}q_1 \\
& \qquad -\frac{\be^2}{2}q_1(1-m_1)\nu \left[  \text{Tanh}^2(g_0, g_1)\right] -   \frac{1}{m_0} \log \EE_0\left[ \exp\left( m_0 f_0(g_0; \boldsymbol q, \boldsymbol m)\right) \right]\,.
\eea \eeq
The computation of the first $\alpha$-derivative is streamlined by the use of Derrida-Ruelle cascades \cite{derrida, ruelle}, as pioneered by Aizenman, Sims and Starr \cite{ass}: consider  a point process 
\beq \eta^{(0)} \defi \left\{ \eta_{i_0}^{(0)}, {i_0 \in \N}\right\}, \eeq 
which is a Poisson-Dirichlet with parameter $m_0$, for short PD$(m_0)$. Independent thereof, and to given $i_0$, consider the point process 
\beq 
\eta^{(1)}_{i_0} \defi \left\{\eta_{i_0, i_1}^{(1)}, {i_1 \in \N}\right\},
\eeq 
which are PD$(m_1)$; these are assumed to be independent over $i_0$, and independent of $\eta^{(0)}$. The Derrida-Ruelle cascade (with two levels) is then the point process
\beq 
\left\{ \eta_{\boldsymbol i}, \boldsymbol i \in \N^2 \right\} \quad \text{where}\quad  \eta_{\boldsymbol i} \defi \eta^{(0)}_{i_0} \eta_{i_0, i_1}^{(1)}\,.
\eeq 
Consider also standard Gaussians $g_{i_0, k}^{(0)}$ and $g_{i_0, i_1, k}^{(1)}$ all independent over $i_0, i_1, k \in \N$. By properties of the Derrida-Ruelle cascades \cite{ass}, the functional \eqref{2map} can be expressed as 
\beq \bea \label{new_repr}
& \Phi_{\alpha}(q_0,  q_1;  \boldsymbol m) =  \\
& =  \frac{1}{N} \E \EE_g \EE_{\eta}\left[ \log \left\{  2^{-N} \sum_{\boldsymbol i, \s \in \Sigma_N} \eta_{\boldsymbol i}  e^{ \sqrt{\alpha} \mathcal H_N(\s) + \sum_{k=1}^N \left( h+ \be \sqrt{q_0} g_{i_0, k}^{(0)}+
\be \sqrt{q_1-q_0} g_{i_0, i_1, k}^{(1)}
 \right) \s_k} \right\} \right],
\eea \eeq
where $\E\EE_g$ stands for expectation of all involved Gaussians, whereas $\EE_\eta$ is expectation of the Derrida-Ruelle cascade. This representation is particularly useful when it comes to the $\alpha$-derivative. In what follows, we will again abuse notations by writing 
\beq
\Phi_\alpha^\star(q_0, q_1; \boldsymbol m) := \Phi_\alpha^\star(q_0^\star(q_0, q_1, \alpha), q_1^\star(q_0, q_1, \alpha); \boldsymbol m)\,,
\eeq
where $q_0^\star, q_1^\star$ are solutions of the corresponding extremality equations, i.e. the $\alpha \neq0$ generalization of \eqref{q0}, \eqref{q1}. It then holds 
\beq \bea \label{onemoretime}
\frac{d}{d\alpha} \Phi_{\alpha}^\star(q_0, q_1 \boldsymbol m) & = - \frac{\partial}{\partial \alpha}  \Phi_{\alpha}(q_0, q_1; \boldsymbol m) \\
& = - \frac{\be^2}{4} + \frac{\be^2}{4} \E \EE_\eta \left< \s_1 \s_2 \tau_1 \tau_2 \right>_{\eta, \alpha}^{\otimes 2} \qquad \text{(Gaussian P.I.)}
\eea \eeq
where $\left< \right>_{\eta, \alpha}^{\otimes 2}$ stands for expectation over the replicated space $
(\N \times \Sigma_N)^2 $ under the Gibbs measure with Hamiltonian \eqref{new_repr} tilted by the Derrida-Ruelle weights. (Remark that in the first line of the above equation, $\boldsymbol q$-extremality plays a key role, insofar it suppresses terms involving $d \boldsymbol q/d\alpha$: these must however be taken into account for the higher order $\alpha$-derivatives). In $\alpha=0$, the spins decouple: again by properties of the Derrida-Ruelle cascades, see e.g. \cite[Chapters 14.2-14.3]{talavol}, one ends up with the neat(er) expression
\beq \bea
\frac{d}{d\alpha} \Phi_{\alpha}^\star(q_0, q_1; \boldsymbol m)\Big|_{\alpha=0} & = 
-\frac{\be^2}{4}+ \frac{\be^2}{4} m_0 \nu\left[ \text{Tanh}(g_0, g_1) \right]^4 \\
& \qquad+ \frac{\be^2}{4}(m_1-m_0)  \nu_0\left[ \nu_1\left[\text{Tanh}(g_0, g_1)\right]^2\right]^2 \\
& \qquad \qquad+ \frac{\be^2}{4}(1-m_1)
\nu\left[\text{Tanh}^2(g_0, g_1)\right]^2 \,.
\eea \eeq
Let us introduce the function 
\beq \bea
G_\text{1rsb}(\boldsymbol q; \boldsymbol m) \defi 
 & \Phi_{0, \boldsymbol m}(q_0, q_1)  -\frac{\be^2}{2}q_0(m_0-m_1) \nu_0 \left[  \nu_1\left[\text{Tanh}(g_0, g_1)\right]^2\right]+ \\
& \qquad \qquad -\frac{\be^2}{2}q_1+ \frac{\be^2}{2}q_1(1-m_1)\nu\left[ \text{Tanh}^2(g_0, g_1)\right] \\
& \qquad +\frac{\be^2}{4}- \frac{\be^2}{4} m_0 \nu\left[\text{Tanh}(g_0, g_1) \right]^4  \\
& \qquad \qquad- \frac{\be^2}{4}(m_1-m_0)  \nu_0\left[ \nu_1\left[ \text{Tanh}(g_0, g_1) \right]^2\right]^2 \\
& \qquad - \frac{\be^2}{4}(1-m_1)
\nu\left[\text{Tanh}^2(g_0, g_1)\right]^2\,,
\eea \eeq
which is the first order expansion of the 1RSB Gibbs potential. By the above,
\beq
F_N(\be, h) \leq \min_{q_0, q_1}  G_\text{1rsb}(\boldsymbol q; \boldsymbol m) + \text{(higher derivatives)}.
\eeq
As in the RS-case, {\it assuming} that the infinite expansions converge absolutely and become irrelevant in the infinite volume limit, we would get
\beq
\lim_N F_N(\be, h) \leq \min_{\boldsymbol q}  G_\text{1rsb}(\boldsymbol q; \boldsymbol m)\,,
\eeq
where the minimum is over $\boldsymbol q$ restricted to the simplex $0 \leq q_0 \leq q_1 \leq 1$.

Remark that we still have freedom in the choice of $\boldsymbol m = (m_0, m_1)$: since the above comes in the form of an upper bound we can thus also minimize w.r.t. the $\boldsymbol m$'s to obtain
\beq \label{1rsb func}
\lim_N F_N(\be, h) \leq \min_{\boldsymbol m} \min_{\boldsymbol q}  G_\text{1rsb}(\boldsymbol q; \boldsymbol m)\,.
\eeq
where $\boldsymbol m, \boldsymbol q$ are both restricted to the respective simplices. 

The above outcome is intriguing: as in the RS-case, the expansion of the Gibbs potential does {\it not} coincide with the Parisi 1RSB functional \cite[p.30]{nishimori}: it does however coincide with the Parisi 1RSB solution if we plug into \eqref{1rsb func} a specific choice of parameters, namely the critical points of the Parisi 1RSB  functional\footnote{These are presumably also critical points of the 1RSB Gibbs potential, but we haven't checked that.}. Precisely:
\begin{itemize}
\item[i)] We choose $q_0, q_1$ to be solution of the fixpoints 
\beq \label{q00}
q_0  = \nu_0\left[ \nu_1[\text{Tanh}(g_0, g_1)]^2\right]
\eeq
and respectively
\beq \label{q11}
q_1= \nu\left[ \text{Tanh}^2(g_0, g_1) \right]\,.
\eeq
(Remark that $0 \leq q_0\leq q_1 \leq 1$ by Jensen's inequality).
\item[ii)] $m_0 :=0$, to be understood in the limiting sense, in which case, as pointed out in \eqref{lim_m} above, the $\nu_0$-density collapses back to the standard Gaussian. 
\end{itemize}
With $m_0=0$,  the parameters $q_0$ and $q_1$ are solutions of the fixed point equations 
\beq
q_0 = \EE_0\left[ \left( \int \tanh(h+ \be \sqrt{q_0} g_0+ \be \sqrt{q_1-q_0} g_1) \nu_1(g_0, dg_1) \right)^2 \right]
\eeq
and respectively, 
\beq
q_1 = \EE_0\left[ \int \tanh^2(h+ \be \sqrt{q_0} g_0+ \be \sqrt{q_1-q_0} g_1) \nu_1(g_0, dg_1) \right]\,.
\eeq
These fixed point equations correspond to the critical points of the 1RSB Parisi functional, cfr. with \cite[Eqs. (3.33) and (3.34)]{nishimori}. Furthermore,  with such choices, the (expansion of the) Gibbs potential takes the form
\beq \bea \label{ultimo}
& G_\text{1rsb}(q_0, q_1; m_0=0, m_1) = \\
& \qquad =  \EE_0\left[ \frac{1}{m_1} \log \EE_1\left[ \cosh\left( h+\be \sqrt{q_0} g_0 + \be \sqrt{q_1-q_0} g_1 \right)^{m_1}\right] \right]  +\frac{\be^2}{2}q_0^2 m_1 + \\
& \qquad \qquad \qquad -\frac{\be^2}{2}q_1+ \frac{\be^2}{2}q_1^2(1-m_1)+\frac{\be^2}{4} - \frac{\be^2}{4}m_1 q_0^2 - \frac{\be^2}{4}(1-m_1)q_1^2 \\
& \qquad = \EE_0\left[ \frac{1}{m_1} \log \EE_1\left[ \cosh\left( h+\be \sqrt{q_0} g_0 + \be \sqrt{q_1-q_0} g_1 \right)^{m_1}\right] \right] + \\
& \hspace{6cm} + \frac{\be^2}{4}\left\{ 1 -2 q_1 + m_1 q_0^2 + q_1^2(1-m_1) \right\}.
\eea \eeq
The minimal value (over $m_1 \in [0,1]$) yields indeed the 1RSB Parisi {\it solution}, cfr. with \cite[Eqn. (3.30)]{nishimori}. We may thus formulate a second bottom line 
\vspace{0.06cm}
\begin{center}
{\sf the 1RSB Parisi solution  is a first-order expansion of the 1RSB-Gibbs potential}\footnote{Thanks to the self-similarity of the building bricks of the Parisi theory an analogous statement is expected for the $K$-RSB solution as well. The formulas  become however so cumbersome that we restrain here from discussing the generic case.}. 
\end{center}
\vspace{0.06cm}
 
In full analogy with the RS-case, one should be aware that what {\it looks} like a first order is again a second order in disguise: in fact, the mean field correction given by the last line of \eqref{ultimo} is again quadratic in $\beta$. This  is again the signature of disorder, which manifests itself through the Gaussian P.I. \eqref{onemoretime}, effectively increasing by one the order of the derivative.

\section{Conclusion and outlook}
This work grew out from our efforts to develop a framework for the analysis of (mean field) spin glasses {\it from first principles}. In the process, we have been strongly influenced by some earlier works of Francesco Guerra, see in particular \cite{guerra_old} where first seeds of an approach to the Parisi theory based on purely thermodynamical considerations have been planted. Prior to this work, Legendre duality principles, albeit in {\it infinite} volume, have  been addressed by Auffinger \& Chen \cite{ac}, who pushed forward a line of research initiated by Guerra in \cite{guerra_l}. It is however unclear to us, what exactly the point of contact between these works and the present notes is: the reason is that in the current setting we perform Legendre transformations on the $\boldsymbol q$-variables to fixed $\boldsymbol m$, whereas in the aforementioned works, the transformations act on both sets of parameters. This issue  definitely deserves further investigations. 

It should be emphasized that the aim here is not that of "yet another proof" of the Parisi formula. For this, a number of rigorous derivations are already available \cite{guerra, tala, panchenko}. Rather, with the pertinent and still prevailing words of Michael Aizenman, the goal is "an analysis which is both rigorous and based on physically recognizable principles" \cite{aizenman}. We hope that the present framework, abiding to the rules of  thermodynamics, brings us one step closer to this goal. 

Of course, the crucial input towards a rigorous implementation of the above program is the convergence of the infinite expansions, be it for the RS- or the RSB-Gibbs potential (or, for that matter, for {\it any} model whatsoever). When/if this is the case we do not know. In fact, given the complexity of the SK-model, one should even brace for the worst, namely that there exist $(\be, h)$-phases where high temperature expansions of the $K$-RSB Legendre functional are always divergent as long as $K< \infty$.
On a positive note one should mention that convergence is only required in a neighborhood of the relevant critical points. But all these issues seem to be premature anyhow:  the much awaited Feynman rules for mean field expansions have been very recently identified by K\"uhn \& Helias in great generality \cite{kh}, but taming the infinite series by means of a rigorous, asymptotical analysis of the expansion-coefficients is a daunting task even for the most "trivial" system, namely the CW-model \cite{kss}. 

Also of  interest is the following issue. In order to build a proper Legendre transformation we have relied on the {\it Ansatz} about the law of the effective fields which stems from the replica computations in the infinite volume limit. It would be very interesting to see how these laws emerge from finite volume. This, of course, is an extremely challenging problem. Some help might come again from the Parisi theory, which is "a solution-facilitating ansatz about the hierarchical form of the replica symmetry breaking" \cite{ass}. It is therefore tempting to believe that the current simplifying Ansatz should be replaced by tree-like Gaussian fields in finite volumes giving rise, in the thermodynamical limit, to the all-important Derrida-Ruelle cascades. Unfortunately, this {\it still} is a challenging problem. In fact, the present framework  acts on finite systems:  we thus need Gaussian fields  which approximate the Derrida-Ruelle cascades in the thermodynamical limit, much akin to the original Derrida's GREM \cite{derrida, derrida_two} and their variants with an infinite number of levels \cite{bovier_kurkova, bovier_kurkova_two}. It is however extremely subtle to construct Gaussian fields in finite volume whose  (approximate) hierarchical structure is encoded by the overlap on Ising configurations. This delicate issue will be addressed elsewhere. 

It should also be mentioned that the framework sketched above potentially lends itself to a $\PP$-almost sure analysis, i.e. with quenched SK-disorder: for this one may consider, say, the 1RSB Legendre formulation \eqref{2map} without performing the $\E$-average. If judiciously combined with the TAP-Plefka treatment \cite{tap, plefka},  this approach shall "enhance" the latter, thereby constructing, and outgoing from finite volumes, TAP-like equations for the ancestors \cite{mpv}. Needless to say, a $\PP$-a.s. analysis would compound the difficulties by orders of magnitude. 

Lastly, we emphasize that concave Legendre transformations such as those implemented above lead to concave envelopes, hence in general to upper bounds. For equalities, one would need a priori knowledge of strict concavity of the finite size Legendre functionals: this is a potentially difficult issue, but there might be a way to bypass the need of such {\it a priori} information. In fact, something mysterious is at work: here we have considered the functional $\boldsymbol q \mapsto \Phi_{\alpha}(\boldsymbol q;  \boldsymbol m)$ and its concave Legendre transformation to {given} $\boldsymbol m$. One can however turn things upside down by considering the functional $\boldsymbol m \mapsto \Phi_{\alpha}(\boldsymbol q;  \boldsymbol m)$, {\it to given $\boldsymbol q$}, and consider the functional which acts on the $\boldsymbol m$-reciprocals, to wit:
\beq 
\boldsymbol t = (t_0,  t_1) \in 1 \leq t_1  \leq t_0 \leq \infty \mapsto \Phi_{\alpha}\left(\boldsymbol q; \frac{1}{t_0}, \frac{1}{t_1} \right).
\eeq 
By Jensen inequality, this functional turns out to be convex (in fact for any $\alpha$, and any $\mathcal H_N$-Hamiltonian). This also implies that the functional $\boldsymbol m \mapsto \Phi_{\alpha}(\boldsymbol q;  \boldsymbol m)$  has a unique minimizer, for any $N$. One can therefore envision an approach via {\it convex} Legendre transformations (be it in the $\boldsymbol t$-, or in the $\boldsymbol m$-formulation, but to fixed $\boldsymbol q$), and this in turns seems to suggest the existence of reversed varational principles, and therefore also of min-max principles, in large but finite volumes. Such min-max principles are known to appear naturally in the cavity fields perturbations of the GREM \cite{bokis_one, bokis_two}; the present take would also elucidate aspects of those (at first sight artificial) models, and their (still opaque) link with the Parisi solution. \\
 
We believe all the above issues are worthy of future research efforts. Here we conclude with a final comment, intentionally cryptic, concerning the title of these notes: the framework worked out above provides evidence to the claim that, after all, 
\begin{center}
{Parisi \emph{is}  Boltzmann}. 
\end{center}

\subsection*{Acknowledgements} It is our pleasure to thank Hermann Dinges for help and advice on convex geometry. NK also wishes to express his enormous gratitude to Timm Plefka for much needed guidance through the conceptual subtleties of thermodynamics, and for sharing his deep insights on mean field spin glasses. Inspiring discussions with Yan V. Fyodorov are also gratefully acknowledged. Last but not least, our thanks to Pierluigi Contucci, Francesco Guerra, Emanuele Mingione, and Dmitry Panchenko for valuable inputs, criticism, and feedback.

\end{document}